\documentclass{amsart}
\usepackage{amsmath}
\usepackage{amssymb}
\usepackage{color}
\usepackage{ifthen}

\newtheorem{theor}{Theorem}
\newtheorem{propo}[theor]{Proposition}
\newtheorem{lemma}[theor]{Lemma}
\newtheorem{corol}[theor]{Corollary}
\theoremstyle{definition}
\newtheorem{defin}[theor]{Definition}
\theoremstyle{remark}
\newtheorem{remar}[theor]{Remark}
\newtheorem{examp}[theor]{Example}

\newtheorem{conve}[theor]{Convention}
\newcommand{\is }{A}

\newcommand{\End }{\mathrm{End}}
\newcommand{\fk }{\mathbf{k}}
\newcommand{\GL }{\mathrm{GL}}
\newcommand{\id }{\mathrm{id}}
\newcommand{\mfg }{\mathfrak{g}}
\newcommand{\mfn }{\mathfrak{n}}
\newcommand{\modu }{\mathrm{mod}\ }
\newcommand{\ndN }{\mathbb{N}}
\newcommand{\ndR }{\mathbb{R}}
\newcommand{\ndZ }{\mathbb{Z}}
\newcommand{\re}{\mathrm{re}}
\newcommand{\rs}{R}
\newcommand{\s}{\sigma }
\newcommand{\sa}{\triangleright}
\newcommand{\sig }[2]{#1\sa #2}
\newcommand{\srs}{\mathcal{R}}
\newcommand{\tC}{\tilde{C}}
\newcommand{\te}{\tilde{e}}
\newcommand{\ts}{\tilde{s}}
\newcommand{\vecm}{\mathbf{m}}
\newcommand{\Wmon}{\widetilde{W}}

\newcommand{\Dchainthree}[2]{
\rule[-6\unitlength]{0pt}{8\unitlength}
\begin{picture}(26,5)(0,3)
\put(1,2){\circle{2}}
\ifthenelse{\equal{#1}{1}}{
 \put(2,2){\line(1,0){10}}}{
  \put(2,2.2){\line(1,0){10}}
  \put(2,1.8){\line(1,0){10}}
  \ifthenelse{\equal{#1}{2l}}{
   \put(1.7,1.1){$\langle $}
  }{
   \put(11,1.1){$\rangle $}
  }
}
\put(13,2){\circle{2}}
\ifthenelse{\equal{#2}{1}}{
 \put(14,2){\line(1,0){10}}}{
  \put(14,2.2){\line(1,0){10}}
  \put(14,1.8){\line(1,0){10}}
  \ifthenelse{\equal{#2}{2l}}{
   \put(13.7,1.1){$\langle $}
  }{
   \put(23,1.1){$\rangle $}
  }
}
\put(25,2){\circle{2}}
\end{picture}}

\newcommand{\Dtriangle}{
\begin{picture}(18,16)(0,3)
\put(4,4){\circle{2}}
\put(5,4){\line(1,0){8}}
\put(14,4){\circle{2}}
\put(4.4472,4.8944){\line(1,2){4.1056}}
\put(9,14){\circle{2}}
\put(13.5528,4.8944){\line(-1,2){4.1056}}
\end{picture}}

\title[Generalization of Coxeter groups]{%
A generalization of Coxeter groups,\\
root systems, and Matsumoto's theorem}

\author{I.~Heckenberger}
\address{Istv\'an Heckenberger, Mathematisches Institut, Universit\"at Leipzig,
         Augustusplatz 10-11, 04109 Leipzig, Germany.}
\email{heckenberger@math.uni-leipzig.de}

\author{H.~Yamane}
\address{
Hiroyuki Yamane,
Department of Pure and Applied Mathematics,
Graduate School of Information Science
and Technology, Osaka University, Toyonaka 560-0043,
Japan
}
\email{yamane@ist.osaka-u.ac.jp}

\begin{document}

\begin{abstract}
 The root systems appearing in the theory of Lie superalgebras and
 Nichols algebras admit a large symmetry extending properly
 the one coming from the Weyl group. Based on this observation we
 set up a general framework in which the symmetry object is a groupoid.
 We prove that in our context the groupoid is generated by reflections and
 Coxeter relations. This answers a question of Serganova.
 Our weak version of the exchange condition allows us
 to prove Matsumoto's theorem. Therefore the word problem is solved for the
 groupoid.
\end{abstract}

\maketitle

\section{Introduction}

Kac-Moody algebras \cite{b-Kac90}
and their generalizations enjoy continuously growing
interest since their introduction in the late sixties.
Lie superalgebras \cite{a-Kac77} became central objects in physical models,
and generalized Kac-Moody-Borcherds algebras \cite{a-Borch88}
turned out to be
important in group theory in connection with the monster.
One of the basic features of these algebras is their relation to
root systems, Weyl groups, and their generalizations.

In connection with contragredient Lie superalgebras Serganova 
\cite{a-Serg96} found out that
the Weyl group symmetry of the Lie algebra can be extended
using reflections on isotropic roots. These so called ``odd reflections''
were used by the second author \cite{a-Yam99} to analyze the structure of
affine super Lie algebras very efficiently.
Serganova \cite{a-Serg96}
gave an axiomatic definition of ``generalized root systems'' (GRS)
which is based on ordinary and odd reflections, and classified them
completely. Below Example~6.7 she writes the following.
\begin{itemize}
\item
``It is an interesting question what is the analogue of Weyl group for
GRS. It is natural to consider the finite group $\Wmon $ generated by
$s_\alpha $ for all $\alpha \in R$. This group in general is not a
subgroup of $\GL (V)$. It is unknown if this group is Coxeter.''
\end{itemize}
Our paper gives an answer to this question as explained in
Remark~\ref{rm:question}.

The discovery of quantized Kac-Moody algebras by Drinfel'd \cite{inp-Drinfeld1}
and Jimbo \cite{a-Jimbo1} opened a new major research direction
with many novel features, ideas, and connections to other research
fields. Motivated by these examples, Andruskiewitsch and Schneider
\cite{a-AndrSchn98}
launched a project to classify pointed Hopf algebras with certain
finiteness conditions, and performed it successfully under some hypotheses
\cite{a-AndrSchn05p}. Their method is based on the knowledge of the
finiteness property of Nichols algebras, the latter being main ingredients for
the construction of quantized Kac-Moody algebras. These finiteness
properties can be determined effectively using the Weyl groupoid
\cite{a-Heck04c} attached to a Nichols algebra. Also the wish to understand
better
the structure of this groupoid gave motivation for our work.

In the literature there exist several generalizations of root systems,
mostly depending on the applications. A very nice general treatment
with many references is the work of Loos and Neher \cite{a-LoosNeh05}.
In our paper we do not attempt to give a complete set of axioms
for a root system. Instead we concentrate on its property being
invariant under symmetry. Motivated by Serganovas work and the
appearence of the Weyl groupoid for Nichols algebras we require
a large symmetry in the sense that reflections with respect to all
simple roots should be defined. Unfortunately this
excludes already from the beginning
the root systems of generalized Kac-Moody-Borcherds algebras.
However Kac-Moody algebras and a very large class of their super and quantized
analogues are still covered, together with so far less understood
examples appearing in the classification of Nichols algebras of
diagonal type \cite{a-Heck06b}. Our main results are Theorems~\ref{th:ffrep},
\ref{th:wexc} and \ref{th:matsu}.
Theorem~\ref{th:ffrep} gives an answer to Serganovas question about the right
generalization of the Weyl group (see also Remark~\ref{rm:question}).
The groupoid
attached naturally to the root system of a Lie superalgebra
is a Coxeter groupoid. The result holds more generally for any root
system satisfying the axioms required in Definition~\ref{df:srs}.
Further, Coxeter groups are known to satisfy the exchange condition.
In our generality this property does not hold in its standard form.
However Theorem~\ref{th:wexc} gives a weak version of it, which seem to
be sufficient in many applications. In particular we are able to
prove Matsumoto's theorem \cite{a-Matsu64}, see Theorem~\ref{th:matsu},
which states that
any two reduced expressions of the same length can be transformed
to each other using Coxeter relations only. This means that if the
groupoid is finitely generated then with finite calculation it is
possible to check whether two given elements in the groupoid are
equal. Hence the word problem is solved for any groupoid appearing
in Theorem~\ref{th:matsu}. Finally, in the last section
of the paper we give an example which appears in the classification of
Nichols algebras but is related neither to a
Kac-Moody nor to a Lie superalgebra.

In this paper let $\ndN $ denote the set of positive integers,
and $\ndN _0=\ndN \cup \{0\}$. The composition sign for general groups
will usually be omitted.

\section{Group actions}

For any set $N$ let $F(N)$ denote the free group generated by
the elements of $N$ and their inverses. Further, let $F_2(N)$ denote
the free group generated by the elements of $N$ as involutions.
With other words, $F_2(N)$ is the quotient of $F(N)$ by the
subgroup consisting of products of elements $gn^2g^{-1}$ and $gn^{-2}g^{-1}$,
where $g\in F(N)$ and $n\in N$. Since in $F_2(N)$ one has $n=n^{-1}$
for all $n\in N$,
$F_2(N)$ is generated as a monoid
by the elements of $N$. The unit of $F_2(N)$ is
denoted by $1$.

Let $N$ and $\is $ be nonempty sets. An action $\sa$ of $F_2(N)$ on $\is $ is
a map $\sa :F_2(N)\times \is \to \is $ such that $1\sa a=a$ and
$g\sa (h\sa a)=(gh)\sa a$ for all $g,h\in F_2(N)$ and $a\in \is $.
By definition of $F_2(N)$ this is equivalent to
\begin{align}\label{eq:sa}
n\sa (n\sa a)=a \qquad \text{for all
$n\in N,a\in \is $.}
\end{align}

An action $\sa :F_2(N)\times \is \to \is $ is called
\textit{transitive}, if for each $a,b\in \is $ there exists
$g\in F_2(N)$ such that $g\sa a=b$.
For an arbitrary action $\sa :F_2(N)\times \is \to \is $
one has a unique decomposition of $\is $ into a disjoint
union of subsets $\is _p$, such that $F_2(N)\sa \is _p\subset \is _p$
and the restriction of $\sa $ to $F_2(N)\times \is _p$ is transitive for
all $p$.

Let $N$ be a nonempty set and $N'\subset N$ a nonempty subset.
Then $F_2(N')$ can be regarded as a subgroup of $F_2(N)$.
Let $\sa $ be an action of $F_2(N)$ on $\is $.
Then $\sa $ induces an action $\sa ':F_2(N')\times \is \to \is $ which
is called the restriction of $\sa $ to $F_2(N')\times \is $.
The transitivity of $\sa $ does not imply transitivity of $\sa '$.

Let $\sa :F_2(N)\times \is \to \is $ be an action and $n,n'\in N$.
For each $a\in \is $ define
\begin{align*}
\Theta (n,n';a):=&\{(nn')^m\sa a,\,
(n'n)^m\sa a\,|\,m\in \ndN _0\}.
\end{align*}

Let $\theta (n,n';a):=|\Theta (n,n';a)|$,
the cardinality of $\Theta (n,n';a)$,
which is either in $\ndN $ or is $\infty $.
One obviously has $\Theta (n,n';a)=\Theta (n',n;a)$
and $\theta (n,n';a)=\theta (n',n;a)$.
Moreover, if $\theta (n,n';a)\in \ndN $
then $\Theta (n,n';a)=\{(nn')^m\sa a\,|\, 0\le m<\theta (n,n';a)\}$.
Further, $\Theta (n,n';n\sa a)=\Theta (n,n';n'\sa a)$
by definition, and equations $\theta (n,n';a)=
\theta (n,n';n\sa a)=\theta (n,n';n'\sa a)$ hold.

Let $a_0:=a$, $b_0:=a$, and define recursively
$a_{m+1}:=n\sa b_m$, $b_{m+1}:=n'\sa a_m$ for all $m\in \ndN _0$.
Using Equation~\eqref{eq:sa}
one obtains that
\begin{align}\label{eq:lsta}
 \theta (n,n';a):=
 \begin{cases}
  \infty & \text{if $a_m\not=b_m$ for all $m\in \ndN $,}\\
  \min \{m\in \ndN \,|\,a_m=b_m\} & \text{otherwise.}
 \end{cases}
\end{align}

\section{Coxeter groupoids}

First we give a generalization of the notion of Coxeter groups to
groupoids. For the notion of groupoids consult for example
\cite[Sect.~3.3]{b-ClifPres61}.

\begin{defin}\label{df:Cox}
Let $N$ and $\is $ be nonempty sets and let
$\sa $ be a transitive action of $F_2(N)$ on $\is $.
For each $a\in \is $ and $i,j\in N$
with $i\not=j$ let $m_{i,j;a}$ be a multiple of $\theta (i,j;a)$
lying in $\ndN \setminus \{1\}\cup \infty $.
Set $\vecm:=(m_{i,j;a}\,|\,i,j\in N,i\not=j,a\in \is )$.
Let $W$ be the
groupoid generated by elements $e_a$ and $s_{i,a}$, where
$a\in \is $ and
$i\in N$, and the following relations:
\begin{gather}\label{eq:erel}
\begin{aligned}
e_a^2=&e_a, & e_ae_b=&0 \text{ for }a\not=b,\\
e_{\sig{i}{a}}s_{i,a}=&s_{i,a}e_a=s_{i,a}, & s_{i,\sig{i}{a}}s_{i,a}=&e_a,
\end{aligned}\\
\label{eq:coxrel}
\begin{aligned}
s_is_j\cdots s_js_{i,a}=s_js_i\cdots s_is_{j,a} 
 \text{ ($m_{i,j;a}$ factors)}
 & \qquad \text{if $m_{i,j;a}$ is finite and odd,}\\
s_js_i\cdots s_js_{i,a}=s_is_j\cdots s_is_{j,a}
 \text{ ($m_{i,j;a}$ factors)}
 & \qquad \text{if $m_{i,j;a}$ is finite and even,}
\end{aligned}
\end{gather}
where $a$ is running over all elements of $\is $.
The quintuple $(W,N,\is ,\sa ,\vecm )$ will be called
a \textit{Coxeter groupoid}.
\end{defin}

Note that in our convention the multiplication of two
elements in a groupoid is always defined, but the product may be
the unique element $0$ which satisfies the property $g0=0g=0$
for all $g$. Moreover,
in \eqref{eq:coxrel}
and later on we use the following convention.

\begin{conve}\label{cv:noindices}
Let $(W,N,\is ,\sa ,\vecm )$ be a Coxeter groupoid,
$(i_1,\ldots ,i_m)$ a
sequence of elements in $N$, and
$(b_1,\ldots ,b_m)$ a sequence of elements in $\is $.
By Relations \eqref{eq:erel} the element
$w=s _{i_1,b_1}\cdots s_{i_{m-1},b_{m-1}}s_{i_m,b_m}\in W$
is obviously zero if one of the relations
$b_j=i_{j+1}\sa b_{j+1}$, where $j\ge 1$, fails.
Thus $w\not=0$ implies that $b_j=i_{j+1}\cdots i_m\sa b_m$
for all $j<m$,
and in the latter case we will write
$s_{i_1}\cdots s_{i_{m-1}}s_{i_m,b_m}$ for $w$.
\end{conve}

Define $\ell :W\to \ndN _0\cup \{-\infty \}$ to be the function
such that $\ell (0)=-\infty $, $\ell (e_a)=0$ for all $a\in \is $,
and 
\begin{align*}
 \ell (w)=\min \{m\in \ndN \,|\,&w=s_{i_1}\cdots s_{i_{m-1}}s_{i_m,a}
 \text{ for some $i_1,\ldots ,i_m\in N$, $a\in \is $}\}
\end{align*}
otherwise. One has
\begin{align}\label{eq:lwinv}
 \ell (w)=&\ell (w^{-1}),\\
\label{eq:lww'}
 \ell (ww')\le &\ell (w)+\ell (w')
\end{align}
for each $w,w'\in W$.
One says that a product $w=s_{i_1}\cdots s_{i_{m-1}}s_{i_m,a}\in W$
is \textit{reduced}, if $m=\ell (w)$.

Note that if the cardinality of $\is $ is 1, then
the action $\sa $ is necessarily trivial, and hence
$\theta (i,j;a)=1$ for $a\in \is $ and all $i,j\in N$.
Thus in this case the definition
of $(W,N,\is ,\sa ,\vecm )$ coincides with the definition of a Coxeter group.

\begin{remar}
Coxeter groups are dealt with effectively using reflections obtained
by conjugation of generators. Note that the structure
constants $m_{i,j;a}$
of a Coxeter groupoid depend also on the elements of $\is $,
and hence in general the adjoint action of $W$ on itself is not defined.
Therefore the standard proofs can not always be
generalized to our setting.
\end{remar}

\begin{remar}
Coxeter groups satisfy the exchange condition, which follows from
the existence of sufficiently many reflections. The exchange condition
in its standard form does not hold for arbitrary Coxeter groupoids,
see Section~\ref{sec:example} for an example.
A weak version of the exchange condition is given in Theorem~\ref{th:wexc}
for a special class of Coxeter groupoids.
\end{remar}

We end this section with the definition of special elements of $W$ and
giving some commutation rules which will be needed later.

Following the notation in Definition~\ref{df:Cox}
assume that $a\in \is $ and $i,j\in N$, where $i\not=j$.
Define
\begin{align}\label{eq:defC}
C_{i,j;a}=\begin{cases}
s_is_js_i\cdots s_is_{j,a}
\text{ ($m_{i,j;a}-1$ factors)} &
\text{if $m_{i,j;a}$ is odd,} \\
s_is_js_i\cdots s_js_{i,a}
\text{ ($m_{i,j;a}-1$ factors)} &
\text{if $m_{i,j;a}$ is even.}
\end{cases}
\end{align}
Then Equation~\eqref{eq:coxrel} implies the following relations.
\begin{align}\label{eq:Crel}
s_iC_{j,i;a}=&\begin{cases}
C_{j,i:j\sa a}s_{j,a} & \text{if $m_{j,i;a}$ is odd,}\\
C_{j,i:i\sa a}s_{i,a} & \text{if $m_{j,i;a}$ is even.}
\end{cases}
\end{align}

\section{A generalization of root systems}

This section is devoted to the definition of a generalization
of root systems which admits the symmetry of a Coxeter groupoid.
The required large symmetry parallels the main idea behind root systems
of semisimple Lie algebras. In the examples given below we will show how
classical objects fit into our definition.

\begin{defin}\label{df:srs}
Let $\srs $ be the set consisting of all
triples $(R,N,\is ,\sa )$ such that the following conditions hold.
\begin{enumerate}
  \item \label{en:transact}
    $N$ and $\is $ are sets and $\sa $ is a transitive action
    of $F_2(N)$ on $\is $.
  \item \label{en:Rpi}
    Let $V_0=\ndR ^{|N|}$. Then $R =\{(R_a,\pi _a,S_a)\,|a\in \is \}$, where
    $\pi _a=\{\alpha _{n,a}\,|\,n\in N\}\subset R_a\subset V_0$,
    and $\pi _a$ is a basis of $V_0$ for all $a\in \is $.
  \item \label{en:R+}
    $\rs _a=\rs ^+_a\cup -\rs ^+_a$ for all $a\in \is $, where
    $\rs ^+_a=\ndN _0\pi _a\cap \rs _a$.
  \item \label{en:rbeta}
    For any $i\in N$ and $a\in \is $ one has
    $\ndR \alpha _{i,a}\cap R_a=\{\alpha _{i,a},-\alpha _{i,a}\}$.
  \item \label{en:Sa}
    $S_a=\{\s_{i,a}\,|\,i\in N\}$, and for each $a\in \is $ and
    $i\in N$ one has
    $\s_{i,a}\in \GL(V_0)$, $\s_{i,a}(\rs _a)=\rs _{\sig{i}{a}}$,
    $\s_{i,a}(\alpha _{i,a})=-\alpha _{i,\sig{i}{a}}$, and
    $\s_{i,a}(\alpha _{j,a})\in
    \alpha _{j,\sig{i}{a}}+\ndN _0\alpha _{i,\sig{i}{a}}$ for all $j\in
    N\setminus \{i\}$.
  \item \label{en:tautau}
    $\s_{i,\sig{i}{a}}\s_{i,a}=\id $ for $a\in \is $ and $i\in N$.
  \item \label{en:taueq}
    Let $a\in \is $, $i,j\in N$, $i\not=j$, and
    $d=|\ndN _0\{\alpha _{i,a},\alpha _{j,a}\}\cap \rs _a|$.
    If $d$ is finite then $\theta (i,j;a)$ is finite and it divides $d$.
\end{enumerate}
The cardinality of $N$ is called the
\textit{rank of} $(R,N,\is ,\sa )\in \srs $.
\end{defin}

For any $(R,N,\is ,\sa )\in \srs $ and any $a\in \is $ define
\begin{equation}
\begin{split}
R_a^\re =\{\s _{i_1}\cdots \s _{i_{m-1}}\s_{i_{m,b}}
(\alpha _{j,b})\,|\,&
m\in \ndN _0,b\in \is ,\\
 & i_1,\ldots ,i_m,j\in N, i_1\cdots i_m\sa b=a\}.
\end{split}
\end{equation}
Similarly to Definition~\ref{df:srs}\eqref{en:R+}
we will use the notation
$(R^\re _a)^+=\ndN _0\pi _a\cap R^\re _a$.
Definition~\ref{df:srs}\eqref{en:tautau} implies that for each
$i\in N$ and $a\in \is $ one has $\s _{i,a}(R_a^\re )=R_{i\sa a}^\re $.
Note that this definition of ``real roots'' coincides with
the standard definition in case of Kac-Moody algebras. In general
however our set $R^\re _a$ contains also
``isotropic roots'' (see examples below). The reason for this difference
is the presence of a larger symmetry and the non-existence of a
bilinear form.

\begin{conve}\label{cv:noindices2}
For products of $\s_{i,a}$ we will use a similar convention as
for products of $s_{i,a}$, see Convention~\ref{cv:noindices}.
\end{conve}

\begin{examp}
Let $(R,N,\is ,\sa )\in \srs $ be of rank one. Then
Definition~\ref{df:srs}\eqref{en:transact}
implies that $\is $ has cardinality
1 or 2. Further, by
Definition~\ref{df:srs}\eqref{en:Rpi},\eqref{en:rbeta} one has
$\pi _a=\{\alpha _a\}$ and $R_a=\{\alpha _a,-\alpha _a\}$
for each $a\in \is $, where $\alpha _a\in \ndR \setminus \{0\}$.
The other axioms do not give additional restrictions, and hence one
has precisely two different elements in $\srs $ which have rank 1.
\end{examp}

\begin{examp}\label{ex:KM}
Let $\is =\{a\}$ and $N=\{1,\ldots ,n\}$ for some $n\in \ndN $.
The set $A$ admits the trivial action $\sa $
of $F_2(N)$ given by $i\sa a:=a$ for all $i\in N$, $a\in \is $.
Let $\mfg $ be a Kac-Moody algebra of rank $n$
and let $(R,\pi )$ and
$(R^\re ,\pi )$ be the corresponding set of roots and real roots,
respectively.
Let $s_{i,a}:=s_i$ denote the
reflection with respect to $\alpha _i\in \pi $, and define
$S=\{s_1,\ldots ,s_n\}$.
Then both $((R,\pi ,S),N,\is ,\sa )$ and
$((R^\re ,\pi ,S),N,\is ,\sa )$ are in $\srs $.
\end{examp}

\begin{examp}\label{ex:super}
Let $\mfg$ be a finite dimensional
contragredient Lie superalgebra of rank $n$
and $\widetilde{R}$ its root system.
Define $R=\widetilde{R}\setminus 2\widetilde{R}$
and $N=\{1,\ldots ,n\}$.
Let $\{\pi_a=\{\alpha _{1,a},\ldots ,\alpha _{n,a}\}\,|\,a\in A\}$ be the
set of all ordered bases of $R$, where $A$ denotes an index set (that is
each basis appears $n!$ times, each time with a different ordering).
For any $i\in N$ and $a\in A$ let $r_{i,a}$ be the reflection on
$\alpha _{i,a}$ if $\alpha _{i,a}$ is an even root, and the odd reflection
if $\alpha _{i,a}$ is an odd root. By \cite[Lemma~6.4]{a-Serg96}
there exists a unique
$b\in A$ such that $r_{i,a}(\pi _a)=\pi _b$,
$\alpha _{i,a}=-\alpha _{i,b}$, and $\alpha _{j,a}\in
\alpha _{j,b}+\ndN _0\alpha _{i,b}$ for all $j\not=i$. Define $i\sa a:=b$
and $\s _{i,a}:=\id $.
Then $\sa $ is an action of $F_2(N)$ on $A$, and all axioms (2)--(7)
in Definition~\ref{df:srs} are fulfilled.
Let $\psi $ denote the map from $A$ to the set of all bases of $R$ defined
by ignoring the order.
Let $A'$ be an orbit of $\sa $ in $\is $. By \cite[Lemma~6.4]{a-Serg96}
$\psi |_{A'}$ is surjective and
and by Lemma~\ref{le:wprop}(iii) it is injective.
By definition of $A'$ one gets
$(R',N,A',\sa )\in \srs $, where $R'=\{(R_a,\pi _a,S_a)\,|\,a\in A'\}$.
Note that the cardinality of $A$
can be reduced further by identifying
$a,b\in A'$ if
$(\alpha_{i,a},\alpha_{j,a})=
(\alpha_{i,b},\alpha_{j,b})$ 
and ${\rm{deg}}(\alpha_{i,a})
={\rm{deg}}(\alpha_{i,b})$ for all $i,j\in \{1,\ldots ,n\}$
(for a somewhat different example see also Section~\ref{sec:example}).
Then Serganovas odd reflections in \cite[(1), (2) of Section 6]{a-Serg96}
are reinterpreted in our setting as base changes in $R$.
\end{examp}

\begin{examp}\label{ex:Nichols}
In \cite{a-Heck04c} and \cite{a-Heck04e} arithmetic root systems are
defined. They fit into the setting of Definition~\ref{df:srs} as sketched
below.
These arithmetic root systems are closely related to Nichols algebras,
and under some additional assumptions to quantized enveloping algebras.
Finite arithmetic root systems are fully classified, see \cite{a-Heck06b}
and the references therein.

Let $\fk $ be a field of characteristic zero, $n\in \ndN $,
$E=\{\alpha _1,\ldots ,\alpha _n\}$ an ordered basis of $\ndZ ^n$, and
$q_{ij}\in \fk \setminus \{0\}$, where $q_{ii}\not=1$,
for all $i,j\in \{1,\ldots ,n\}$.
Let $\chi :\ndZ ^n\times \ndZ ^n\to \fk \setminus \{0\}$ be the unique
bicharacter such that $\chi (\alpha _i,\alpha _j)=q_{ij}$
for all $i,j\in \{1,\ldots ,n\}$.
Assume now that for any $i\in
\{1,\ldots ,n\}$ and any $j\not=i$ relation
$(q_{ii}^m-1)(q_{ii}^mq_{ij}q_{ji}-1)=0$ holds for some $m\in \ndN _0$,
and denote the smallest such number by $m_{ij}$.
Let $i\sa E$ denote the ordered basis of $\ndZ ^n$ consisting of
$i\sa \alpha _i=-\alpha _i$ and
$i\sa \alpha _j:=\alpha _i+m_{ij}\alpha _j$ for $j\not=i$,
in the ordering induced by the one of $E$.
Proceed inductively with the new bases, using now the values of the
bicharacter $\chi $ instead of the $q_{ij}$.
The assumption that the bases $i_1\sa (i_2 \sa \cdots (i_m \sa E))$
always exist is crucial to obtain the setting in Definition~\ref{df:srs},
and no general algebraic condition on the $q_{ij}$
is known which is equivalent to this property.
Let $\is $ denote the set of all ordered bases of $\ndZ ^n$ of the form
$i_1\sa (i_2\sa \cdots (i_m\sa E))$.
By definition, $\sa $ is a transitive action of $F_2(\{1,\ldots ,n\})$
on $\is $.

The Nichols algebra attached to the matrix $(q_{ij})$ is $\ndZ ^n$-graded,
where the generators have degree $\alpha _i$, and $i\in \{1,\ldots ,n\}$.
Moreover it has a restricted
PBW basis \cite{a-Khar99}
consisting of homogeneous elements, and the degrees of these PBW
generators can be regarded as the set of positive roots $\Delta ^+$
with respect to the basis $E$.
In particular, let $C=(c_{ij})_{i,j=1,\ldots ,n}$
be a symmetrizable Cartan matrix,
and $d_i\in \ndN$ for $1\le i\le n$ such that $(d_ic_{ij})$ is symmetric.
Let $q\in \fk $ be not a root of $1$, and set $q_{ij}:=q^{d_ic_{ij}}$.
Then the Nichols algebra attached to the structure constants $q_{ij}$
is exactly $U_q(\mfn _+)$, where $\mfn _+$ is the upper
triangular part of the Kac-Moody algebra associated to $C$,
and the set $\Delta ^+$ is exactly the set of positive roots associated
to $C$.

Let $\Delta :=\Delta ^+\cup -\Delta ^+$,
$R_a:=\Delta $ and $s_{i,a}:=\id $ for all $a\in \is $ (see above) and
$i\in \{1,\ldots ,n\}$.
Then $s_{i,a}$ satisfies the conditions in
Definition~\ref{df:srs}\eqref{en:Sa}. The remaining axioms follow
either from the construction or
from the theory of Nichols algebras, in particular from
\cite[Proposition\,1]{a-Heck04c}.

In case of $U_q(\mfn _+)$ it is sufficient to consider a one-element set
$\is $, but with our definition of $\is $ this is not the case.
However, in general it is possible to reduce the number of elements
of $\is $ by introducing
the following equivalence relation. Call
$a=\{\beta _1,\ldots ,\beta _n\}$ and $b=\{\gamma _1,\ldots ,\gamma _n\}$
equivalent, if
$\chi (\beta _i,\beta _i)=\chi (\gamma _i,\gamma _i)$ and
$\chi (\beta _i,\beta _j)\chi (\beta _j,\beta _i)=
\chi (\gamma _j,\gamma _i)\chi (\gamma _i,\gamma _j)$ for all
$i,j\in \{1,\ldots ,n\}$.
Since the matrix of $s_{i,a}$ with respect to the bases
$a$ and $i\sa a$ depends only on the values of $\chi $,
it depends only on the equivalence class of $a$. A similar property
holds for $R_a$, and hence one may replace $\is $ by the set of its equivalence
classes. If $q_{ii}$ is not a root of $1$ for all $i$
(in particular in the above setting related to $U_q(\mfn _+)$)
then with this identification
one obtains exactly the setting in Example~\ref{ex:KM}. A very concrete
example of different type can be found in Section~\ref{sec:example}.
\end{examp}

\begin{remar}
In view of the large amount of examples given above
it seems very difficult to determine all elements of $\srs $.
On the other hand, Theorem~\ref{th:matsu} shows that the structure of the
elements of $\srs $ is not as complicated as it looks at first sight.
We see the rich internal structure together with the difficulty of
classification as the fascinating key features of our
definition.
\end{remar}

We continue with analyzing the structure of the elements of $\srs $.

\begin{lemma}\label{le:sigperm}
For any $a\in \is $ and $i\in N$ one has
$\s _{i,a}(R^+_a\setminus \{\alpha _{i,a}\})
=R^+_{i\sa a}\setminus \{\alpha _{i,i\sa a}\}$.
\end{lemma}

\begin{proof}
 This follows from Definition~\ref{df:srs}\eqref{en:R+},%
 \eqref{en:rbeta},\eqref{en:Sa}.
\end{proof}

\begin{corol}
  Let $m\in \ndN $, $(i_1,\ldots ,i_m)\in N^m$, and $a\in \is $. Let
  $b=i_1\cdots i_m\sa a$ and $w:=\s _{i_2}\cdots \s _{i_{m-1}}\s _{i_m,a}$.
  Then one has
  \begin{align*}
    &\s _{i_1,i_1\sa b}w(R^+_a)\cap -R^+_b=\\
    &\quad \begin{cases}
      \s _{i_1,i_1\sa b}(w(R^+_a)\cap -R^+_{i_1\sa b})\cup
      \{-\alpha _{i_1,b}\} & \text{if $\alpha _{i_1,i_1\sa b}\in w(R^+_a)$}\\
      \s _{i_1,i_1\sa b}( (w(R^+_a)\cap -R^+_{i_1\sa b})\setminus
      \{-\alpha _{i_1,i_1\sa b}\}) & \text{otherwise,
      that is $-\alpha _{i_1,i_1\sa b}\in w(R^+_a)$.}
    \end{cases}
  \end{align*}
  \label{co:swRcap-R}
\end{corol}

\begin{corol}\label{co:postoneg}
Let $m\in \ndN _0$, $(i_1,\ldots ,i_m)\in N^m$, and $a\in \is $.
For all $r\le m$
define $a_r:=i_r\cdots i_m \sa a$, and
let $w:=\s _{i_1}\cdots \s _{i_{m-1}}\s _{i_m,a}$.
Then one has $|w(R^+_a)\cap -R^+_{a_1}|\in m-2\ndN _0$.
Moreover if $|w(R^+_a)\cap -R^+_{a_1}|=m$ then
\begin{align}\label{eq:wR+a}
w(R^+_a)\cap -R^+_{a_1}=\{-\s _{i_1}\cdots \s_{i_{r-2}}
\s _{i_{r-1},a_r}(\alpha _{i_r,a_r})\,|\,1\le r\le m\}.
\end{align}
\end{corol}

Let $N$ and $\is $ be nonempty sets and $\sa $ an action of
$F_2(N)$ on $\is $.
Let $N'$ be a nonempty subset of $N$.
Then there is a unique decomposition
$\is =\cup _p\is _p$ of $\is $
into the disjoint union of nonempty subsets such that
$F_2(N')$ acts transitively on $\is _p$ for all $p$. Let
$\sa _p$ denote the restriction of $\sa $ to
$F_2(N')\times \is _p$.

The following lemma can be obtained immediately from
Definition~\ref{df:srs}.

\begin{lemma}\label{le:restr}
Let $(R=\{(R_a,\pi _a,S_a)\,|\,a\in \is \},N,\is ,\sa )\in \srs $
and $N'$ a nonempty subset of $N$. Let
$\is =\cup _p\is _p$ be the above decomposition with respect to
the action of $F_2(N')$.
For any $a\in \is $ set
\begin{align}
\pi '_a:=&\{\alpha _{n,a}\,|\,n\in N'\},&
S'_a:=&\{\s _{n,a}\,|\,n\in N'\},&
R'_{p,a}:=&\ndZ \pi '_a\cap R_a.
\end{align}
Then $(R'_p,N',\is _p,\sa _p)\in \srs $, where
$R'_p=\{(R'_{p,a},\pi '_a,S'_a)\,|\,a\in \is _p\}$.
\end{lemma}



\section{The rank two case}

In this section let $(R,N=\{i,j\},\is ,\sa )\in \srs $
have rank two. We will give a description of $R^+_a$,
where $a\in \is $,
in terms of elements of some $\pi _b$, $b\in \is $.
Set $d:=|R^+_a|\in \ndN \cup \{\infty \}$.

\begin{lemma} \label{le:no--}
 Let $a\in \is $, $m\in \ndN _0$,
 $(i_1,\ldots ,i_m)\in N^m$, and
 $w=\s _{i_1}\cdots \s _{i_{m-1}}\s _{i_m,a}$.
 If $m<d$ and $w(\alpha _{i,a})\in -R^+_b$,
 where $b=i_1\cdots i_m\sa a$,
 then $w(\alpha _{j,a})\in R^+_b$.
\end{lemma}

\begin{proof}
 Assume to the contrary that $w(\alpha _{j,a})\in -R^+_b$.
 Then $w(R^+_a)\subset -R^+_b$ which is a contradiction to
 Corollary~\ref{co:postoneg} and the relation $m<d$.
\end{proof}
%
%

\begin{lemma}\label{le:R+a}
 Let $a\in \is $. For all $n\in \ndN _0$
 set $i_{2n}=j$, $i_{2n+1}=i$,
 $j_{2n}=i$, $j_{2n+1}=j$, $a_n=i_n\cdots i_1\sa a$, and
$b_n=j_n\cdots j_1\sa a$.

(i) If $d\in \ndN $ then
 \begin{align}\label{eq:R+a}
  R^+_a=\{\s _{i_1}\cdots \s _{i_{m-1}}\s _{i_m,a_m}(\alpha _{i_{m+1},a_m})
  \,|\, 0\le m<d\}.
 \end{align}
In this case one has $\alpha _{j,a}=\s _{i_1}\cdots \s _{i_{d-2}}\s _{i_{d-1},
 a_{d-1}}(\alpha _{i_d,a_{d-1}})$.

(ii) If $R_a$ is infinite then
 \begin{align}\label{eq:infR+a}
 (R^\re _a)^+=
 \{
 \s _{i_1}\cdots \s _{i_{m-1}}\s _{i_m,a_m}(\alpha _{i_{m+1},a_m})\,|\,
 m\in \ndN _0\}
 \cup {}&\\ \notag
 \{
 \s _{j_1}\cdots \s _{j_{m-1}}\s _{j_m,b_m}(\alpha _{j_{m+1},b_m})\,|\,
 & m\in \ndN _0\},
 \end{align}
and the elements given in Equation~\eqref{eq:infR+a} are pairwise
different from each other.
\end{lemma}

\begin{proof}
Set $w_0:=\id $, $w_1:=\s _{i,a}$, $M_0:=\emptyset $, and
$M_1:=\{\alpha _{i,a}\}$.
We inductively define
$w_t\in \End (\ndR ^2)$ and $M_t\subset R^+_a$
for $t\ge 0$ by the formulas
\begin{align}
 w_{t+1}:=&\s _{i_{t+1},a_t}w_t,&
 M_{t+1}:=&M_t\cup \{w_t^{-1}(\alpha _{i_{t+1},a_t})\}.
\end{align}
We will show that
\begin{enumerate}
\item \label{en:Mt}
 $M_t=\{\beta \in R^+_a\,|\,w_t(\beta )\in -R^+_{a_t}\}$,
\item
 \label{en:cardMt}
 $|M_t|=t$.
\end{enumerate}
Setting $t=d$ in the above formulas and using the inductive
definition of $M_d$ one obtains Equation~\eqref{eq:R+a}.
Similarly,
Equation~\eqref{eq:infR+a} follows from the definition of
$R^\re _a$
and Definition~\ref{df:srs}(6).
The remaining statements
of (i) and (ii) follow from the fact that
$\alpha _{j,a}\notin M_t$ for $t<d$, which
itself is a consequence of
the assumption $\alpha _{i,a}\in M_t$ and Lemma~\ref{le:no--}.

We perform the proof of \eqref{en:Mt}
and \eqref{en:cardMt} by induction
on $t$.
Note that for $t=0$ and $t=1$ these formulas hold trivially.
Assume now that $M_t$ satisfies \eqref{en:Mt} and \eqref{en:cardMt},
and that relation $t<d$ holds.
Then $w_t^{-1}(\alpha _{i_t,a_t})=w_{t-1}^{-1}\s _{i_t,a_{t-1}}^{-1}
(\alpha _{i_t,a_t})
=-w_{t-1}^{-1}(\alpha _{i_t,a_{t-1}})\in -M_{t-1}$, so
$w_t^{-1}(\alpha _{i_t,a_t})\in -R^+_a$.
Hence $w_t^{-1}(\alpha _{i_{t+1},a_t})\in R^+_a$ by Lemma~\ref{le:no--}.
This gives $\alpha _{i_{t+1},a_t}\in w_t(R^+_a)$.
By Corollary~\ref{co:swRcap-R} one has
\begin{equation} \label{eq:t->t+1}
  \begin{aligned}
  R^+_a\cap -w_{t+1}^{-1}(R^+_{a_{t+1}})=&w^{-1}_{t+1}(w_{t+1}(R^+_a)
  \cap -R^+_{a_{t+1}})\\
  =&w^{-1}_{t+1}(\s _{i_{t+1},a_t}(w_t(R^+_a)\cap -R^+_{a_t})\cup 
  \{-\alpha _{i_{t+1},a_{t+1}}\})\\
  =&w^{-1}_t(w_t(R^+_a)\cap -R^+_{a_t})\cup 
  \{-w^{-1}_{t+1}(\alpha _{i_{t+1},a_{t+1}})\}\\
  =&M_t\cup \{w^{-1}_t(\alpha _{i_{t+1},a_t})\}=M_{t+1}.
  \end{aligned}
\end{equation}
Hence \eqref{en:Mt} holds. Moreover, by the second equality of
\eqref{eq:t->t+1} and Corollary~\ref{co:swRcap-R} one can see that
\eqref{en:cardMt} is valid.
\end{proof}

\begin{lemma}\label{le:rank2cox}
If $d$ is finite then for all $a\in \is $ the relations
 \begin{equation}\label{eq:rank2cox}
 \begin{aligned}
 \s_i\s_j\s_i\cdots \s_j\s_{i,a}=
 \s_j\s_i\s_j\cdots \s_i\s_{j,a}
 \text{ ($d$ factors)}
 & \quad \text{if $d$ is odd,}\\
 \s_j\s_i\s_j\cdots \s_j\s_{i,a}=
 \s_i\s_j\s_i\cdots \s_i\s_{j,a}
 \text{ ($d$ factors)}
 & \quad \text{if $d$ is even}
 \end{aligned}
 \end{equation}
hold.
\end{lemma}

\begin{proof}
 By Definition~\ref{df:srs}\eqref{en:tautau},\eqref{en:taueq}
 the assertion of the lemma is equivalent to the equation
 $(\s_i\s_j)^{d-1}\s_i\s_{j,a}=\id $. If $d$ is even then
 applying Lemma~\ref{le:R+a}(i) twice gives that
 \begin{align*}
  &(\s_i\s_j)^{d-1}\s_i\s_{j,a}(\alpha _{i,a})=
  (\s_i\s_j)^{d/2}\s_{i,i \sa a}(\alpha _{i,i \sa a})\\
  &\qquad =(\s_i\s_j)^{d/2-1}\s_i\s_{j,a}(-\alpha _{i,a})
  =\s _{i,i \sa a}(-\alpha _{i,i \sa a})=\alpha _{i,a},\\
  &(\s_i\s_j)^{d-1}\s_i\s_{j,a}(\alpha _{j,a})=
  -(\s_i\s_j)^{d-1}\s_{i,j \sa a}(\alpha _{j,j \sa a})\\
  &\qquad =-(\s_i\s_j)^{d/2-1}\s_i\s_{j,a}(\alpha _{j,a})
  =(\s_i\s_j)^{d/2-1}\s_{i,j \sa a}(\alpha _{j,j\sa a})
  =\alpha _{j,a}.
 \end{align*}
 For odd $d$ the claim follows similarly.
\end{proof}



\section{The general case}

In this section we use the results in the rank two case to
show that to any $(R,N,\is ,\sa )\in \srs $ one can attach
a faithful representation of a Coxeter groupoid in a natural way.

\begin{propo}\label{pr:rep}
Let $(R,N,\is ,\sa )\in \srs $. For each $i,j\in N$,
where $i\not=j$, and any $a\in \is $
set
$m_{i,j;a}:=|\ndN _0\{\alpha _{i,a},\alpha _{j,a}\}\cap \rs _a|$.
Set $V=\bigoplus _{a\in \is }V_a$, where $V_a=V_0=\ndR ^{|N|}$ (see
Definition~\ref{df:srs}\eqref{en:Rpi}), and let
$P_a:V\to V_a$ and $\iota _a:V_a\to V$
be the canonical projection
and canonical injection, respectively. Then the assignment
$\rho :e_a\mapsto \iota _aP_a$,
$s_{i,a}\mapsto \iota _{\sig{i}{a}}\s _{i,a}P_a$,
gives a representation $(\rho ,V)$ of the Coxeter groupoid
$(W,N,\is ,\sa ,\vecm )$.
In particular for $W$ one has $e_a\not=0$ and $s_{i,a}\not=0$
for all $i\in N$ and $a\in \is $.
\end{propo}

\begin{proof}
By construction and by Definition~\ref{df:srs}\eqref{en:tautau}
one obtains that Equations~\eqref{eq:erel} are compatible
with the definition of the representation.
By Definition~\ref{df:srs}\eqref{en:taueq} it remains to show
that if $m_{i,j;a}$ is finite then equation
\begin{align}\label{eq:speccox}
\begin{aligned}
\s_i\s_j\s_i\cdots \s_j\s_{i,a}=
\s_j\s_i\s_j\cdots \s_i\s_{j,a}
 \text{ ($m_{i,j:a}$ factors)}
 & \quad \text{if $m_{i,j;a}$ is odd,}\\
\s_j\s_i\s_j\cdots \s_j\s_{i,a}=
\s_i\s_j\s_i\cdots \s_i\s_{j,a}
 \text{ ($m_{i,j;a}$ factors)}
 & \quad \text{if $m_{i,j:a}$ is even}
\end{aligned}
\end{align}
holds. The restrictions of the above equations to the space
$\ndR \alpha _{i,a}\oplus \ndR \alpha _{j,a}$ are valid
by Lemma~\ref{le:rank2cox} and Lemma~\ref{le:restr}.
This means that for $x_{i,j;a}:=
(\s_i\s_j)^{m_{i,j;a}-1}\s_i\s_{j,a}$ one has
\begin{align*}
x_{i,j;a}(\alpha _{i,a})=\alpha _{i,a},\quad
x_{i,j;a}(\alpha _{j,a})=\alpha _{j,a},\quad
x_{i,j;a}(\alpha _{n,a})\in \alpha _{n,a}
 +\ndN_0\{\alpha _{i,a},\alpha _{j,a}\}
\end{align*}
for all $n\in N\setminus \{i,j\}$.
Since $x_{i,j;a}(R_a)=R_a$
by Definition~\ref{df:srs}\eqref{en:Sa},\eqref{en:taueq},
this implies that $x_{i,j;a}(\alpha _{n,a})=\alpha _{n,a}$
for all $n\in N$.
\end{proof}

If $(R,N,\is ,\sa )\in \srs $ then let
$(\rho (W),R,N,\is ,\sa )$ denote the groupoid obtained as
the image of the representation $\rho $ in Proposition~\ref{pr:rep}.


In the following lemmata let
 $(R,N,\is ,\sa )\in \srs $, and for all $a\in \is $ and      
 $i,j\in N$, where $i\not=j$, set                                 
 $m_{i,j;a}:=|\ndN _0\{\alpha _{i,a},\alpha _{j,a}\}              
 \cap \rs _a|$.
Note that one has
\begin{align}
m_{i,j;a}=m_{j,i;a}=m_{i,j;i\sa a}=m_{i,j;j\sa a}.
\end{align}
We will analyze the properties of the Coxeter groupoid
$(W,N,\is ,\sa ,\vecm )$.

\begin{lemma}
 Let $m\in \ndN _0$, $(i_1,\ldots ,i_m)\in N^m$, and $a\in \is $.
 If $\ell (s_{i_1}\cdots s_{i_{m-1}}s_{i_m,a})<m$
 then $s_{i_2}\cdots s_{i_{m-1}}s_{i_m,a}=
 s_{i_1}s_{j_2}\cdots s_{j_{m-2}}s_{j_{m-1},a}$
 for some $j_2,\ldots ,j_{m-1}\in N$.
\end{lemma}

\begin{proof}
By the definition of $\ell $, the relation
$\ell (s_{i_1}\cdots s_{i_{m-1}}s_{i_m,a})<m$
implies that equation $s_{i_1}\cdots s_{i_{m-1}}s_{i_m,a}=
s_{j_2}\cdots s_{j_{r-1}}s_{j_r,a}$ holds for some $r\le m$.
Moreover $W$ is $\ndZ /2\ndZ $-graded, and hence
$m$ is even if and only if $r$ is odd. Using
Equations~\eqref{eq:erel} one can assume that $r=m-1$.
Then the claim follows from the last equation in \eqref{eq:erel}.
\end{proof}

\begin{lemma}\label{le:wprop}
Let $m\in \ndN $, $(i_1,\ldots ,i_m)\in N^m$, and $a\in \is $.

 (i) If $\s _{i_1}\cdots \s _{i_{m-2}}\s _{i_{m-1},a}
 (\alpha _{i_m,a})\in -R^+_b$, where $b=i_1\cdots i_{m-1}\sa a$,
 then relation
 $\ell (s_{i_1}\cdots s_{i_{m-1}}s_{i_m,i_m\sa a})<m$ holds.

 (ii) If $\ell (s_{i_1}\cdots s_{i_{m-1}}s_{i_m,a})=m$,
 and there exist $k\in N\setminus \{i_1\}$,
 $\alpha \in \pi _a$, and $h_1,h_2\in \ndN $ such that
 $\s_{i_2}\cdots \s_{i_{m-1}}\s_{i_m,a}(\alpha )=
 h_1\alpha _{i_1,f}+h_2\alpha _{k,f}$, where
 $f=i_2\cdots i_m \sa a$, then equation
 $s_{i_2}\cdots s_{i_{m-1}}s_{i_m,a}=
 s_ks_{j_2}\cdots s_{j_{m-2}}s_{j_{m-1},a}$
 holds for some $j_2,\ldots ,j_{m-1}\in N$.

 (iii) One has $\ell (s_{i_1}\cdots s_{i_{m-1}}s_{i_m,a})
 =|\s _{i_1}\cdots \s_{i_{m-1}}\s_{i_m,a}(R^+_a)
 \cap -R^+_{i_1\cdots i_m \sa a}|$.
\end{lemma}

\begin{proof}
We prove all statements parallelly by induction on $m$.
Note that for $m=1$ the claim holds by Definition~\ref{df:srs}.
Assume now that the lemma holds for $0,1,\ldots ,m-1$.

To (i): By Equation~\eqref{eq:lww'} one can assume that
$\ell (s_{i_2}\cdots s_{i_{m-1}}s_{i_m,i_m \sa a})=m-1$.
Moreover the last equation in \eqref{eq:erel} yields that
it is sufficient to consider the case $i_1\not=i_2$.
By the converse of (i) for $m-1$ one has
$\s_{i_2}\cdots \s_{i_{m-2}}\s_{i_{m-1},a}(\alpha _{i_m,a})\in
R^+_{i_1 \sa b}$.
By Lemma~\ref{le:sigperm}, since $\s _{i_1,i_1\sa b}(\alpha _{i_1,i_1\sa b})
=-\alpha _{i_1,b}$, one has
\begin{align}\label{eq:walpha}
\s_{i_2}\cdots \s_{i_{m-2}}\s_{i_{m-1},a}(\alpha _{i_m,a})
=\alpha _{i_1,i_1 \sa b}.
\end{align}
We prove now by induction, with letting $d:=m_{i_1,i_2;b}\in (2+\ndN _0)\cup
\{\infty \}$, that 
\begin{align} \label{eq:replacei}
  s_{i_1}s_{i_2}\cdots s_{i_{m-1}}s_{i_m,i_m\sa a}=
  s_{i_1}s_{i_2}s_{i_3 (\modu 2)}\cdots s_{i_d (\modu 2)}s_{j_{d+1}}
  \cdots s_{j_{m-1}}s_{i_m,i_m\sa a}
\end{align}
for some $j_{d+1}\cdots j_{m-1}\in N$, where $n (\modu 2)\in \{1,2\}$
for $n\le d$. We also show that $m>d$,
and that for these $j_{d+1},\ldots ,j_{m-1}$ relation
\begin{align}\label{eq:ssalpha}
\s _{j_{d+1}}\cdots
\s _{j_{m-2}}\s _{j_{m-1},a}(\alpha _{i_m,a})=
\alpha _{i_{d+1 (\modu 2)},a'},
\end{align}
holds, where $a'=j_{d+1}\cdots j_{m-1}\sa a$.
%
Thus in our setting we may apply Equation~\eqref{eq:coxrel} to get
$$s_{i_1}\cdots
s_{i_{m-1}}s_{i_m,i_m\sa a}=s_{i_2}s_{i_1}s_{i_2}\cdots s_{i_{d+1 (\modu 2)}}
s_{j_{d+1}}\cdots s_{j_{m-1}}s_{i_m,i_m\sa a}.$$
Since $\s_{i_{d+1 (\modu 2)}}\s_{j_{d+1}}\cdots \s_{j_{m-2}}\s_{j_{m-1},a}
(\alpha _{i_m,a})=\s_{i_{d+1 (\modu 2),a'}}(\alpha _{i_{d+1 (\modu 2)},a'})$
is contained in $-R^+_{i_{d+1 (\modu 2)}\sa a'}$ and because of $d\ge 2$
we may apply (i) to conclude that
$\ell (s_{i_1}\cdots s_{i_{m-1}}s_{i_m,i_m\sa a})<m$.

Note that equation~$i_n=i_{n (\modu 2)}$ holds trivially for $n\le 2$.
Assume now that it also holds for
$n\le p$, where $2\le p<d$.
Then Equation~\eqref{eq:walpha} and Lemma~\ref{le:R+a}
imply that
$$
\s _{i_{p+1}}\cdots \s_{i_{m-2}}\s_{i_{m-1},a}(\alpha _{i_m,a})
=\underset{\text{$p-1$ factors}}
{\underbrace{\s _{i_{p (\modu 2)}}\s _{i_{p-1 (\modu 2)}}
\cdots \s _{i_1} \s _{i_2,i_1\sa b}}}(\alpha _{i_1,i_1 \sa b})
$$
is an element in $R^+_{i_{p+1}\cdots i_{m-1}\sa a}\cap
\ndN _0\{\alpha _{i_1,i_{p+1}\cdots i_{m-1}\sa a},
\alpha _{i_2,i_{p+1}\cdots i_{m-1}\sa a}\}$, but it is not in
$\pi _{i_{p+1}\cdots i_{m-1}\sa a}$. In particular $m-1\ge p+1$.
Therefore induction hypothesis (ii)
can be applied to the element $s_{i_p}s_{i_{p+1}}
\cdots s_{i_{m-2}}s_{i_{m-1},a}\in W$ and taking
$k:=i_{p+1 (\modu 2)}$ and $\alpha :=\alpha _{i_m,a}$.
One obtains that $s_{i_{p+1}}\cdots s_{i_{m-2}}s_{i_{m-1},a}
=s_{i_{p+1} (\modu 2)}s_{j'_{p+2}}\cdots s_{j'_{m-2}}s_{j'_{m-1},a}$
for some $j'_{p+2},\ldots ,j'_{m-1}\in N$. Thus induction on $p$ gives that
Equation~\eqref{eq:replacei} holds and that $m-1\ge d$.

Equation~\eqref{eq:ssalpha} follows from Equation~\eqref{eq:walpha}
and Lemma~\ref{le:R+a} using \eqref{eq:replacei}.

To (iii):
If $\ell (s_{i_1}\cdots s_{i_{m-1}}s_{i_m,a})<m$ then
the claim follows from the induction hypothesis (iii) and
Proposition~\ref{pr:rep}. So one may assume that
$\ell (s_{i_1}\cdots s_{i_{m-1}}s_{i_m,a})=m$ and hence
$\ell (s_{i_2}\cdots s_{i_{m-1}}s_{i_m,a})=m-1$.
Corollary~\ref{co:postoneg} gives that
\begin{align*}
 |\s _{i_1}\cdots \s_{i_{m-1}}\s_{i_m,a}(R^+_a)
 \cap -R^+_{i_1\cdots i_m \sa a}|\in m-2\ndN _0.
\end{align*}
Assume that (iii) does not hold. The induction
hypothesis (iii) and Corollary~\ref{co:swRcap-R} imply that
$-\alpha _{i_1,i_1\sa c}=\s _{i_2}\cdots \s_{i_{m-1}}\s_{i_m,a}(\beta )$
for some $\beta \in R^+_a$, where $c=i_1\cdots i_m\sa a$.
Hence
$\s _{i_m}\cdots \s_{i_3}\s_{i_2,i_1 \sa c}
(\alpha _{i_1,i_1 \sa c})\in -R^+_a$.
Since (i) is already proven for $m$, we conclude that
$\ell (s_{i_m}\cdots s_{i_2}s_{i_1,c})<m$ which is
a contradiction to the assumption and Equation~\eqref{eq:lwinv}.

To (ii):
By (iii) and Corollary~\ref{co:postoneg} one has
$-\alpha _{i_1,i_1 \sa f}=
\s_{i_1}\cdots \s_{i_{m-1}}\s_{i_m,a}(\beta )$
for some $\beta \in R^+_a$, and hence
$\s_{i_2}\cdots \s_{i_{m-1}}\s_{i_m,a}(\beta )=
\alpha _{i_1,f}$.
Suppose that
$\beta ':=
\s _{i_m}\cdots \s _{i_3}\s _{i_2,f}(\alpha _{k,f})\in R^+_a$.
Then
$\s_{i_2}\cdots \s_{i_{m-1}}\s_{i_m,a}(\beta ')=
\alpha _{k,f}$, and hence
$$\s_{i_2}\cdots \s_{i_{m-1}}\s_{i_m,a}(h_1\beta +h_2\beta ')=
h_1\alpha _{i_1,f}+h_2\alpha _{k,f}=
\s_{i_2}\cdots \s_{i_{m-1}}\s_{i_m,a}(\alpha ).$$
Since $h_1\beta +h_2\beta '\notin \pi _a$ but
$\alpha \in \pi _a$, this gives a contradiction.
Therefore relation
$\s_{i_m}\cdots \s_{i_3}\s_{i_2,f}(\alpha _{k,f})\in -R^+_a$ 
holds, and hence $\ell (s_{i_m}\cdots s_{i_2}s_{k,k \sa f})<m$ by (i).
Equation~\eqref{eq:lwinv} and the $\ndZ /2\ndZ $-grading
of $W$ give that $s_ks_{i_2}\cdots s_{i_{m-1}}s_{i_m,a}=
s_{j_2}\cdots s_{j_{m-2}}s_{j_{m-1},a}$ for some
$j_2,\ldots ,j_{m-1}\in N$. Multiplying the last equation
with $s_k$ from the left
gives the claim.
\end{proof}

\begin{corol}\label{co:l=m-1}
Let $m\in \ndN $, $(i_1,\ldots ,i_m,j)\in N^{m+1}$,
and $a\in \is $, and suppose that
$\ell (s_{i_1}\cdots s_{i_{m-1}}s_{i_m,a})=m$. Then
$\ell (s_{i_1}\cdots s_{i_m}s_{j,j\sa a})=m-1$ if and only if
$\s_{i_1}\cdots \s_{i_{m-1}}\s_{i_m,a}(\alpha _{j,a})\in
-R^+_{i_1\cdots i_m\sa a}$. Equivalently, one has
$\ell (s_{i_1}\cdots s_{i_m}s_{j,j\sa a})=m+1$ if and only if
$\s_{i_1}\cdots \s_{i_{m-1}}\s_{i_m,a}(\alpha _{j,a})\in
R^+_{i_1\cdots i_m\sa a}$.
\end{corol}

\begin{proof}
  Equation~\eqref{eq:lwinv}, Corollary~\ref{co:swRcap-R}, and 
  Lemma~\ref{le:wprop}(iii) give that
  $\ell (s_{i_1}\cdots s_{i_m}s_{j,j\sa a}) \in \{m-1,m+1\}$.
  Therefore the if part follows from
  Lemma~\ref{le:wprop}(i). On the other hand, if
  $\ell (s_{i_1}\cdots s_{i_m}s_{j,j\sa a})=m-1$
  then
  $\ell (s_js_{i_m}\cdots s_{i_2}s_{i_1,b})=m-1$,
  where $b=i_1\cdots i_m\sa a$.
  Using the assumption
  $\ell (s_{i_m}\cdots s_{i_2}s_{i_1,b})=m$,
  Corollary~\ref{co:swRcap-R} and Lemma~\ref{le:wprop}(iii)
  give that
  $\s_{i_m}\cdots \s_{i_2}\s_{i_1,b}(\beta )= -\alpha _{j,a}$
  for some $\beta \in R^+_b$. This implies the only if part of
  the claim.
\end{proof}

\begin{lemma}\label{le:coxbeg}
Let $m\in \ndN $, $(i_0,i_1,\ldots ,i_m,j)\in N^{m+2}$,
and $a\in \is $, and suppose that
$\ell (s_{i_1}\cdots s_{i_{m-1}}s_{i_m,a})=m$.
If $\s _{i_1}\cdots \s _{i_{m-1}}\s _{i_m,a}(\alpha _{j,a})=
\alpha _{i_0,i_1\cdots i_m\sa a}$ then $i_0\not=i_1$ and
\begin{gather}
s_{i_1}\cdots s_{i_{m-1}}s_{i_m,a}=
s_{i_{1 (\modu 2)}}s_{i_{2 (\modu 2)}}\cdots
s_{i_{d-1 (\modu 2)}}s_{j_d}\cdots s_{j_{m-1}}s_{j_m,a},\\
\s_{j_d}\cdots \s_{j_{m-1}}\s_{j_m,a}(\alpha _{j,a})=
\alpha _{i_{d (\modu 2)},j_d\cdots j_m\sa a}
\end{gather}
for some $j_d,\cdots ,j_m\in N$,
where $d=m_{i_1,i_2;i_1\cdots i_m\sa a}$ and
$k (\modu 2)\in \{0,1\}$ for $k\in \{1,\ldots ,d\}$.
\end{lemma}

\begin{proof}
By Corollary~\ref{co:l=m-1} one has
$\ell (s_{i_1}\cdots s_{i_m}s_{j,j\sa a})=m+1$.
By Lemma~\ref{le:wprop}(iii) and Corollary~\ref{co:postoneg},
especially Equation~\eqref{eq:wR+a},
one has $-\s _{i_1}\cdots \s_{i_{m-1}}\s _{i_m,a}(\alpha _{j,a})\not=
-\alpha _{i_1,i_1\cdots i_m\sa a}$, so $i_1\not=i_0$.
Since $\s _{i_0}\s _{i_1}\cdots \s _{i_{m-1}}\s _{i_m}(\alpha _{j,a})
\in -R^+_{i_0\cdots i_m\sa a}$, one can
apply the proof of Lemma~\ref{le:wprop}(i), especially
Equations~\eqref{eq:walpha}, \eqref{eq:replacei}, and \eqref{eq:ssalpha}.
\end{proof}

\begin{theor}\label{th:ffrep}
 Let $(R,N,\is ,\sa )\in \srs $, and for all $a\in \is $ and      
 $i,j\in N$, where $i\not=j$, set                                 
 $m_{i,j;a}:=|\ndN _0\{\alpha _{i,a},\alpha _{j,a}\}              
 \cap \rs _a|$.
 Then the representation $(\rho ,V)$ of the
 Coxeter groupoid $(W,N,\is ,\sa ,\vecm )$ is faithful.
\end{theor}

\begin{proof}
 Let $a\in \is $, $m\in \ndN _0$, and $(i_1,\ldots ,i_m)\in N^m$.
 One has to show that if $\s _{i_1}\cdots \s_{i_{m-1}}\s _{i_m,a}
 =\id $ and $i_1\cdots i_m\sa a=a$ then $s_{i_1}\cdots s_{i_{m-1}}
 s_{i_m,a}=e_a$. This follows from Lemma~\ref{le:wprop}(iii).
\end{proof}

\begin{remar}\label{rm:question}
Theorem~\ref{th:ffrep} gives one possible answer to Serganova's question
about the analog of Weyl group mentioned in the introduction.
On the one hand we do not consider the group proposed by Serganova.
On the other hand Theorem~\ref{th:ffrep} shows that with our definition
one obtains a Coxeter groupoid attached to any contragredient Lie superalgebra,
which acts faithfully on the root system, and contains all features of
ordinary and odd reflections.
\end{remar}

\begin{corol}\label{co:longestword}
Let $(R,N,\is ,\sa )\in \srs $ and assume that $d:=|R_a|$ is finite for one
(that is each) $a\in \is $. Then $\is $ is finite. Further,
let $(W,N,\is ,\sa ,\vecm )$ be the associated Coxeter groupoid.
Then for each $a\in \is $ there exists a unique
$w_a\in W$
such that $w_ae_a=w_a$ and $\ell (w_a)$ is maximal.
Moreover, one has $\ell (w_a)=d$.
\end{corol}

\begin{proof}
If $R_a$ is finite then $|R_a|$ is independent of $a\in \is $ by
Definition~\ref{df:srs}\eqref{en:Sa}.

Corollary~\ref{co:l=m-1}
implies that if $\ell (w_a)$ is maximal then $w_a(\alpha _{j,a})\in
-R^+_b$ for all $j\in N$, where $b\in \is $ is such that
$e_bw_a=w_a$.
Assume now that $w_a$ and $w'_a$ have maximal length. Then
$w_a(w'_a)^{-1}$ maps $R^+_c$ to $R^+_b$, where $b,c\in A$ are such that
$e_bw_a=w_a$ and $e_cw'_a=w'_a$. By Lemma~\ref{le:wprop}(iii) one gets
$\ell (w_a(w'_a)^{-1})=0$, so $b=c$ and $w_a(w'_a)^{-1}=e_b$.
Hence $w_a=w'_a$.
The existence of a $w_a$ with maximal length follows
also from Lemma~\ref{le:wprop}(iii).
Now the finiteness of $\is $ follows from the transitivity of $\sa $
and the fact that
there are only finitely many elements $w\in W$ with $we_a=w$ and $\ell (w)
\le |R^+_a|$.
\end{proof}

\begin{theor}[weak exchange condition]\label{th:wexc}
Let $m\in \ndN $, $(i_0,i_1,\ldots ,i_m,j)\in N^{m+2}$,
and $a\in \is $, and suppose that
$\ell (s_{i_1}\cdots s_{i_{m-1}}s_{i_m,a})=m$.
If $\s _{i_1}\cdots \s _{i_{m-1}}\s _{i_m,a}(\alpha _{j,a})=
\alpha _{i_0,i_1\cdots i_m\sa a}$ then
there exist $r\in \ndN $, $(j_1,\ldots ,j_r)\in N^r$,
$(k_1,\ldots ,k_{r+1})\in N^{r+1}$, and
$(a_1,\ldots ,a_r)\in \is ^r$, such that the following relations hold.
\begin{align}
\label{eq:normalform}
 s_{i_1}\cdots s_{i_{m-1}}s_{i_m,a}=&
 C_{j_1,k_1;a_1}\cdots C_{j_r,k_r;a_r},\\
\label{eq:wexc1}
 C_{j_1,k_1;a_1}\cdots C_{j_r,k_r;a_r}s_{j,j\sa a}=&
 s_{i_0}C_{j_1,k_1;k_2\sa a_1}\cdots C_{j_r,k_r;k_{r+1}\sa a_r},\\
\label{eq:wexc2}
 C_{j_1,k_1;k_2\sa a_1}\cdots C_{j_r,k_r;k_{r+1}\sa a_r}s_{j,a}=&
 s_{i_0}C_{j_1,k_1;a_1}\cdots C_{j_r,k_r;a_r},
\end{align}
\begin{gather}\label{eq:cons1}
j_1=i_1,\ k_1=i_0,\ k_{r+1}=j,\ a_r=a,\
\sum _{t=1}^rm_{j_t,k_t;a_t}-r=m,
\\
\label{eq:cons2}
j_n\not=k_n\text{ for $1\le n\le r$,}\\
\label{eq:cons3}
 k_{n+1}=
 \begin{cases}
  j_n & \text{if $m_{j_n,k_n;a_n}$ is odd,} \\
  k_n & \text{if $m_{j_n,k_n;a_n}$ is even,}
 \end{cases}
\text{ for $1\le n\le r$,}
\end{gather}
and for $2\le n\le r$ one has
\begin{gather}
\label{eq:cons4}
 a_{n-1}=
 \begin{cases}
  j_nk_n\cdots j_nk_n\sa a_n
  \text{ ($m_{j_n,k_n;a_n}-1$ factors)} &
  \text{if $m_{j_n,k_n;a_n}$ is odd,} \\
  j_nk_n\cdots k_nj_n\sa a_n
  \text{ ($m_{j_n,k_n;a_n}-1$ factors)} &
  \text{if $m_{j_n,k_n;a_n}$ is even.}
 \end{cases}
\end{gather}
\end{theor}

\begin{proof}
We show Equations~\eqref{eq:normalform},\eqref{eq:wexc1}, and
\eqref{eq:wexc2} in a way that \eqref{eq:cons1},
\eqref{eq:cons2}, \eqref{eq:cons3}, and \eqref{eq:cons4} also hold.

Equation~\eqref{eq:normalform} is obtained easily by
iterated application of Lemma~\ref{le:coxbeg}.
By Equations~\eqref{eq:Crel} and \eqref{eq:cons1}-\eqref{eq:cons4},
for $1\le n\le r$, one has
$C_{j_n,k_n;a_n}s_{k_{n+1},k_{n+1}\sa a_n}=s_{k_n}C_{j_n,k_n;k_{n+1}\sa a_n}$.
Hence Equation~\eqref{eq:wexc1} holds.
Equation~\eqref{eq:wexc2} can be obtained from
\eqref{eq:wexc1} by multiplication with $s_{j,a}$ from the right
and $s_{i_0}$ from the left.
\end{proof}

Now we are going to prove a generalization of Matsumoto's theorem
\cite{a-Matsu64}.
First we need a definition.

\begin{defin}
Let $(W,N,\is ,\sa ,\vecm )$ be a Coxeter groupoid.
Let $\Wmon $ denote the semigroup generated by the set
$\{0,\te_a,\ts_{i,a}\,|\,a\in \is ,i\in N\}$
and the relations
\begin{gather}
00=0,\quad 0\te_a=\te_a0=0\ts_{i,a}=\ts_{i,a}0=0,\\
\te_a^2=\te_a,\quad  \te_a\te_b=0 \text{ for }a\not=b,\quad
\te_{\sig{i}{a}}\ts_{i,a}=\ts_{i,a}\te_a=\ts_{i,a},\\
\begin{aligned}
\ts_i\ts_j\cdots \ts_j\ts_{i,a}=\ts_j\ts_i\cdots \ts_i\ts_{j,a} 
 \text{ ($m_{i,j;a}$ factors)}
 & \qquad \text{if $m_{i,j;a}$ is finite and odd,}\\
\ts_j\ts_i\cdots \ts_j\ts_{i,a}=\ts_i\ts_j\cdots \ts_i\ts_{j,a}
 \text{ ($m_{i,j;a}$ factors)}
 & \qquad \text{if $m_{i,j;a}$ is finite and even,}
\end{aligned}
\end{gather}
where a similar convention as in Convention~\ref{cv:noindices}
is used. More precisely we will write
$(\Wmon ,N,\is ,\sa ,\vecm )$ for this semigroup.
\end{defin}

Analogously to Equation~\eqref{eq:defC}
we define
\begin{align}\label{eq:deftC}
\tC _{i,j;a}=\begin{cases}
\ts_i\ts_j\ts_i\cdots \ts_i\ts_{j,a}
\text{ ($m_{i,j;a}-1$ factors)} &
\text{if $m_{i,j;a}$ is odd,} \\
\ts_i\ts_j\ts_i\cdots \ts_j\ts_{i,a}
\text{ ($m_{i,j;a}-1$ factors)} &
\text{if $m_{i,j;a}$ is even.}
\end{cases}
\end{align}
Similarly to Equations~\eqref{eq:Crel},\eqref{eq:wexc1} and \eqref{eq:wexc2}
we have
\begin{gather} \label{eq:Crelt}
  \ts_i\tC_{j,i;a}=\begin{cases}
  \tC_{j,i:j\sa a}\ts_{j,a} & \text{if $m_{j,i;a}$ is odd,}\\
  \tC_{j,i:i\sa a}\ts_{i,a} & \text{if $m_{j,i;a}$ is even,}
  \end{cases}\\
\label{eq:wexc1t}
 \tC_{j_1,k_1;a_1}\cdots \tC_{j_r,k_r;a_r}\ts_{j,j\sa a}=
 \ts_{i_0}\tC_{j_1,k_1;k_2\sa a_1}\cdots \tC_{j_r,k_r;k_{r+1}\sa a_r},\\
\label{eq:wexc2t}
 \tC_{j_1,k_1;k_2\sa a_1}\cdots \tC_{j_r,k_r;k_{r+1}\sa a_r}\ts_{j,a}=
 \ts_{i_0}\tC_{j_1,k_1;a_1}\cdots \tC_{j_r,k_r;a_r},
\end{gather}
where $i,j,j_n,k_n,a_n$ are as in Theorem~\ref{th:wexc}.

\begin{theor}\label{th:matsu}
 Let $(R,N,\is ,\sa )\in \srs $, and for all $a\in \is $ and      
 $i,j\in N$, where $i\not=j$, set                                 
 $m_{i,j;a}:=|\ndN _0\{\alpha _{i,a},\alpha _{j,a}\}              
 \cap \rs _a|$. Suppose that $m\in \ndN _0$, $a\in \is $,
 and $(i_1,\ldots ,i_m),(j_1,\ldots ,j_m)\in N^m$ such that
 $\ell (s_{i_1}\cdots s_{i_{m-1}}s_{i_m,a})=m$ and equation
 $s_{i_1}\cdots s_{i_{m-1}}s_{i_m,a}=
 s_{j_1}\cdots s_{j_{m-1}}s_{j_m,a}$ holds
 in $(W,N,\is ,\sa ,\vecm )$.
 Then in the semigroup
 $(\Wmon ,N,\is ,\sa ,\vecm )$
 one has
 $\ts_{i_1}\cdots \ts_{i_{m-1}}\ts_{i_m,a}=
 \ts_{j_1}\cdots \ts_{j_{m-1}}\ts_{j_m,a}$.
\end{theor}

\begin{proof}
 We proceed by induction on $m$, where for $m=0$ and $m=1$
 the claim holds by Proposition~\ref{pr:rep} and
 Lemma~\ref{le:sigperm}.

 Assume now that the claim holds for some $m\in \ndN $.
 Since $s_{i_1}\cdots s_{i_m}s_{j_m,j_m\sa a}=
 s_{j_1}\cdots s_{j_{m-2}}s_{j_{m-1},j_m\sa a}$, we have
 $\ell (s_{i_1}\cdots s_{i_m}s_{j_m,j_m\sa a})=m-1$.
 By Corollary~\ref{co:l=m-1} this implies that
 $\s _{i_1}\cdots \s _{i_{m-1}}\s _{i_m,a}(\alpha _{j_m,a})
 \in -R^+_{i_1\cdots i_m\sa a}$. Therefore there exists
 $t\in \{1,\dots ,m\}$ such that relations
 $\s _{i_{t+1}}\cdots \s _{i_{m-1}}\s _{i_m,a}(\alpha _{j_m,a})
 \in R^+_{i_{t+1}\cdots i_m\sa a}$ and
 $\s _{i_t}\s _{i_{t+1}}\cdots \s _{i_{m-1}}\s _{i_m,a}
 (\alpha _{j_m,a})\in -R^+_{i_ti_{t+1}\cdots i_m\sa a}$ hold.
 Thus by Lemma~\ref{le:sigperm} one has
 $\s _{i_{t+1}}\cdots \s _{i_{m-1}}\s _{i_m,a}(\alpha _{j_m,a})
 =\alpha _{i_t,i_{t+1}\cdots i_m\sa a}$.
 Now we apply Theorem~\ref{th:wexc} and obtain that there exists
 $r\in \ndN $ and $(h_1,\ldots ,h_r)\in N^r$ such that
 \begin{align}\label{eq:s=C}
  s_{i_{t+1}}\cdots s_{i_{m-1}}s_{i_m,a}=
  C_{h_1,k_1;a_1}\cdots C_{h_r,k_r;a_r}
 \end{align}
 where the same notation as in Theorem~\ref{th:wexc} is used.
 Note that $k_1=i_t$, $h_1=i_{t+1}$, $k_{r+1}=j_m$, and $a_r=a$.
 By Equation~\eqref{eq:wexc1} we obtain that
 \begin{align}\label{eq:th1}
  s_{i_1}\cdots s_{i_{m-1}}s_{i_m,a}s_{j_m,j_m\sa a}=
  s_{i_1}\cdots s_{i_{t-1}}C_{h_1,k_1;k_2\sa a_1}\cdots
  C_{h_r,k_r;k_{r+1}\sa a_r}.
 \end{align}
 By induction hypothesis this yields
 \begin{align}\label{eq:th2}
  \ts_{j_1}\cdots \ts_{j_{m-2}}\ts_{j_{m-1},j_m\sa a}=
  \ts_{i_1}\cdots \ts_{i_{t-1}}\tC_{h_1,k_1;k_2\sa a_1}\cdots
  \tC_{h_r,k_r;k_{r+1}\sa a_r}.
 \end{align}
 Since $t\ge 1$, induction hypothesis and Equation~\eqref{eq:s=C}
 give that
 \begin{align}\label{eq:th3}
  \ts_{i_{t+1}}\cdots \ts_{i_{m-1}}\ts_{i_m,a}=
  \tC_{h_1,k_1;a_1}\cdots \tC_{h_r,k_r;a_r}
 \end{align}
 By Equation~\eqref{eq:wexc2t} this implies that
 \begin{align}\label{eq:th4}
  \ts_{i_t}\cdots \ts_{i_{m-1}}\ts_{i_m,a}=
  \tC_{h_1,k_1;k_2\sa a_1}\cdots
  \tC_{h_r,k_r;k_{r+1}\sa a_r}\ts _{j_m,a}.
 \end{align}
 Now multiply Equation~\eqref{eq:th2} from the right by
 $\ts _{j_m,a}$ and
 Equation~\eqref{eq:th4} from the left by $\ts _{i_1}\cdots
 \ts _{i_{t-1}}$ to obtain the claim.
\end{proof}

\section{An example}
\label{sec:example}

We demonstrate the structural results in the previous sections
on an example related to Nichols algebras, see
\cite[Table 2,row $15$]{a-Heck05a}.

Let $\is =\{a,b,c,d,e\}$ and $N=\{1,2,3\}$. Further, let
$\pi _x:=\{\alpha _{1,x},\alpha _{2,x},\alpha _{3,x}\}$
be a basis of $\ndZ ^3\subset \ndR ^3$
for each $x\in \is $. Define
\begin{align*}
R^+_a:=\{\alpha _{1,a}, \alpha _{2,a},& \alpha _{3,a},
\alpha _{1,a}+\alpha _{2,a}, \alpha _{2,a}+\alpha _{3,a},
2\alpha _{2,a}+\alpha _{3,a}, \alpha _{1,a}+\alpha _{2,a}+\alpha _{3,a},\\
& \alpha _{1,a}+2\alpha _{2,a}+\alpha _{3,a},
\alpha _{1,a}+2\alpha _{2,a}+2\alpha _{3,a},
\alpha _{1,a}+3\alpha _{2,a}+2\alpha _{3,a}\},\\
R^+_b:=\{\alpha _{1,b}, \alpha _{2,b},& \alpha _{3,b},
\alpha _{1,b}+\alpha _{2,b}, \alpha _{2,b}+\alpha _{3,b},
\alpha _{1,b}+2\alpha _{2,b}, 2\alpha _{2,b}+\alpha _{3,b},\\
& \alpha _{1,b}+\alpha _{2,b}+\alpha _{3,b},
\alpha _{1,b}+2\alpha _{2,b}+\alpha _{3,b},
\alpha _{1,b}+3\alpha _{2,b}+\alpha _{3,b}\},\\
R^+_c:=\{\alpha _{1,c}, \alpha _{2,c},& \alpha _{3,c},
\alpha _{1,c}+\alpha _{2,c}, \alpha _{2,c}+\alpha _{3,c},
\alpha _{1,c}+\alpha _{2,c}+\alpha _{3,c},
\alpha _{1,c}+2\alpha _{2,c}+\alpha _{3,c},\\
& 2\alpha _{1,c}+2\alpha _{2,c}+\alpha _{3,c},
\alpha _{1,c}+2\alpha _{2,c}+2\alpha _{3,c},
2\alpha _{1,c}+3\alpha _{2,c}+2\alpha _{3,c}\},\\
R^+_d:=\{\alpha _{1,d}, \alpha _{2,d},& \alpha _{3,d},
\alpha _{1,d}+\alpha _{2,d}, \alpha _{2,d}+\alpha _{3,d},
\alpha _{1,d}+2\alpha _{2,d}, \alpha _{1,d}+\alpha _{2,d}+\alpha _{3,d},\\
& \alpha _{1,d}+2\alpha _{2,d}+\alpha _{3,d},
2\alpha _{1,d}+2\alpha _{2,d}+\alpha _{3,d},
2\alpha _{1,d}+3\alpha _{2,d}+\alpha _{3,d}\},\\
R^+_e:=\{\alpha _{1,e}, \alpha _{2,e},& \alpha _{3,e},
\alpha _{1,e}+\alpha _{2,e}, \alpha _{1,e}+\alpha _{3,e},
\alpha _{2,e}+\alpha _{3,e}, \alpha _{1,e}+\alpha _{2,e}+\alpha _{3,e},\\
& 2\alpha _{1,e}+\alpha _{2,e}+\alpha _{3,e},
\alpha _{1,e}+\alpha _{2,e}+2\alpha _{3,e},
2\alpha _{1,e}+\alpha _{2,e}+2\alpha _{3,e}\}.
\end{align*}
According to the following figure
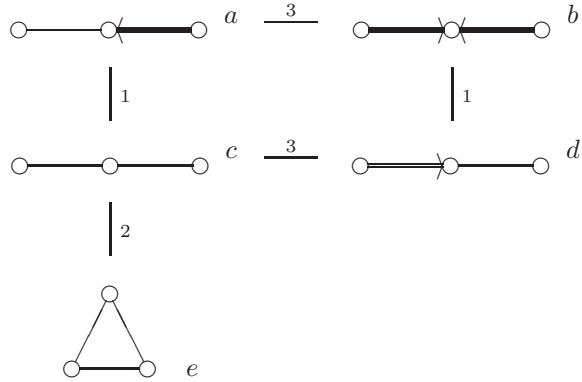
\begin{center}
\setlength{\unitlength}{1mm}
\begin{tabular}{ccc}
\Dchainthree{1}{2l} \ {$a$} & $\overset{3}{\rule{2em}{.5pt}}$ &
\Dchainthree{2r}{2l} \ $b$ \\
\rule[-.5em]{.5pt}{2em}\ $\overset{1}{}$ & &
\rule[-.5em]{.5pt}{2em}\ $\overset{1}{}$\\
\Dchainthree{1}{1} \ $c$ & $\overset{3}{\rule{2em}{.5pt}}$ &
\Dchainthree{2r}{1} \ $d$ \\
\rule[-.5em]{.5pt}{2em}\ $\overset{2}{}$ & & \\
\Dtriangle \ $e$ 
\end{tabular}
\end{center}
define the maps $\s _{i,x}$ for all $i\in \{1,2,3\}$ and $x\in \is $.
Edges between Dynkin diagrams $x,y\in \is $ labelled by the number $i$
mean that $i\sa x=y$ (and hence $i\sa y=x$). If there is no edge with label
$i$ which connects a given Dynkin diagram $x$ with another diagram,
then one has $i\sa x=x$. The numbers $m\in \ndN _0$ in the elements
$\s_{i,x}(\alpha _{j,x})=\alpha _{j,y}+m\alpha _{i,y}$ are defined
in the conventional way. Thus one has
\begin{align*}
\s_{1,a}(\alpha _{1,a})=&-\alpha _{1,c}, &
\s_{1,a}(\alpha _{2,a})=&\alpha _{2,c}+\alpha _{1,c}, &
\s_{1,a}(\alpha _{3,a})=&\alpha _{3,c},\\
\s_{2,a}(\alpha _{2,a})=&-\alpha _{2,a}, &
\s_{2,a}(\alpha _{1,a})=&\alpha _{1,a}+\alpha _{2,a}, &
\s_{2,a}(\alpha _{3,a})=&\alpha _{3,a}+2\alpha _{2,a},\\
\s_{3,a}(\alpha _{3,a})=&-\alpha _{3,b}, &
\s_{3,a}(\alpha _{1,a})=&\alpha _{1,b}, &
\s_{3,a}(\alpha _{2,a})=&\alpha _{2,b}+\alpha _{3,b},
\end{align*}
\begin{align*}
\s_{1,b}(\alpha _{1,b})=&-\alpha _{1,d}, &
\s_{1,b}(\alpha _{2,b})=&\alpha _{2,d}+\alpha _{1,d}, &
\s_{1,b}(\alpha _{3,b})=&\alpha _{3,d},\\
\s_{2,b}(\alpha _{2,b})=&-\alpha _{2,b}, &
\s_{2,b}(\alpha _{1,b})=&\alpha _{1,b}+2\alpha _{2,b}, &
\s_{2,b}(\alpha _{3,b})=&\alpha _{3,b}+2\alpha _{2,b},\\
\s_{3,b}(\alpha _{3,b})=&-\alpha _{3,a}, &
\s_{3,b}(\alpha _{1,b})=&\alpha _{1,a}, &
\s_{3,b}(\alpha _{2,b})=&\alpha _{2,a}+\alpha _{3,a},
\end{align*}
\begin{align*}
\s_{1,c}(\alpha _{1,c})=&-\alpha _{1,a}, &
\s_{1,c}(\alpha _{2,c})=&\alpha _{2,a}+\alpha _{1,a}, &
\s_{1,c}(\alpha _{3,c})=&\alpha _{3,a},\\
\s_{2,c}(\alpha _{2,c})=&-\alpha _{2,e}, &
\s_{2,c}(\alpha _{1,c})=&\alpha _{1,e}+\alpha _{2,e}, &
\s_{2,c}(\alpha _{3,c})=&\alpha _{3,e}+\alpha _{2,e},\\
\s_{3,c}(\alpha _{3,c})=&-\alpha _{3,d}, &
\s_{3,c}(\alpha _{1,c})=&\alpha _{1,d}, &
\s_{3,c}(\alpha _{2,c})=&\alpha _{2,d}+\alpha _{3,d},
\end{align*}
\begin{align*}
\s_{1,d}(\alpha _{1,d})=&-\alpha _{1,b}, &
\s_{1,d}(\alpha _{2,d})=&\alpha _{2,b}+\alpha _{1,b}, &
\s_{1,d}(\alpha _{3,d})=&\alpha _{3,b},\\
\s_{2,d}(\alpha _{2,d})=&-\alpha _{2,d}, &
\s_{2,d}(\alpha _{1,d})=&\alpha _{1,d}+2\alpha _{2,d}, &
\s_{2,d}(\alpha _{3,d})=&\alpha _{3,d}+\alpha _{2,d},\\
\s_{3,d}(\alpha _{3,d})=&-\alpha _{3,c}, &
\s_{3,d}(\alpha _{1,d})=&\alpha _{1,c}, &
\s_{3,d}(\alpha _{2,d})=&\alpha _{2,c}+\alpha _{3,c},
\end{align*}
\begin{align*}
\s_{1,e}(\alpha _{1,e})=&-\alpha _{1,e}, &
\s_{1,e}(\alpha _{2,e})=&\alpha _{2,e}+\alpha _{1,e}, &
\s_{1,e}(\alpha _{3,e})=&\alpha _{3,e}+\alpha _{1,e},\\
\s_{2,e}(\alpha _{2,e})=&-\alpha _{2,c}, &
\s_{2,e}(\alpha _{1,e})=&\alpha _{1,c}+\alpha _{2,c}, &
\s_{2,e}(\alpha _{3,e})=&\alpha _{3,c}+\alpha _{2,c},\\
\s_{3,e}(\alpha _{3,e})=&-\alpha _{3,e}, &
\s_{3,e}(\alpha _{1,e})=&\alpha _{1,e}+\alpha _{3,e}, &
\s_{3,e}(\alpha _{2,e})=&\alpha _{2,e}+\alpha _{3,e}.
\end{align*}
One can check that $\s _i(R_x)=R_{i\sa x}$ for all $i\in \{1,2,3\}$
and $x\in \{a,b,c,d,e\}$.
Proposition~\ref{pr:rep} tells that the following relations hold
in $(W,N,\is ,\sa ,\vecm )$.
\begin{align*}
s_1s_2s_{1,a}=&s_2s_1s_{2,a},&
s_1s_{3,a}=&s_3s_{1,a},&
s_2s_3s_2s_{3,a}=&s_3s_2s_3s_{2,a},\\
s_1s_2s_1s_{2,b}=&s_2s_1s_2s_{1,b},&
s_1s_{3,b}=&s_3s_{1,b},&
s_2s_3s_2s_{3,b}=&s_3s_2s_3s_{2,b},\\
s_1s_2s_{1,c}=&s_2s_1s_{2,c},&
s_1s_{3,c}=&s_3s_{1,c},&
s_2s_3s_{2,c}=&s_3s_2s_{3,c},\\
s_1s_2s_1s_{2,d}=&s_2s_1s_2s_{1,d},&
s_1s_{3,d}=&s_3s_{1,d},&
s_2s_3s_{2,d}=&s_3s_2s_{3,d},\\
s_1s_2s_{1,e}=&s_2s_1s_{2,e},&
s_1s_3s_{1,e}=&s_3s_1s_{3,e},&
s_2s_3s_{2,e}=&s_3s_2s_{3,e}.
\end{align*}

Consider now the following longest word in $W$.
\begin{align*}
&s_{1,b}s_{2,b}s_{1,d}s_{3,c}s_{2,e}
s_{3,e}s_{1,e}s_{3,e}s_{2,c}s_{1,a}=\\
& s_{1,b}s_{2,b}s_{1,d}s_{2,d}s_{3,c}
s_{2,e}s_{1,e}s_{3,e}s_{2,c}s_{1,a}=\\
& s_{2,d}s_{1,b}s_{2,b}s_{1,d}s_{3,c}
s_{2,e}s_{1,e}s_{3,e}s_{2,c}s_{1,a}=\\
& s_{2,d}s_{1,b}s_{2,b}s_{3,a}s_{1,c}
s_{2,e}s_{1,e}s_{3,e}s_{2,c}s_{1,a}=\\
& s_{2,d}s_{1,b}s_{2,b}s_{3,a}s_{2,a}
s_{1,c}s_{2,e}s_{3,e}s_{2,c}s_{1,a}=\\
& s_{2,d}s_{1,b}s_{2,b}s_{3,a}s_{2,a}
s_{1,c}s_{3,d}s_{2,d}s_{3,c}s_{1,a}=\\
& s_{2,d}s_{1,b}s_{2,b}s_{3,a}s_{2,a}
s_{1,c}s_{3,d}s_{2,d}s_{1,b}s_{3,a}.
\end{align*}
The second and the last line are an example for the weak exchange
condition. However, one has the relation
\begin{align*}
&s_{1,b}s_{2,b}s_{1,d}s_{3,c}s_{2,e}
s_{3,e}s_{1,e}s_{3,e}s_{2,c}s_{1,a}\not=
s_{2,c}s_{1,a}s_{3,b}s_{2,b}
s_{3,a}s_{1,c}s_{3,d}s_{2,d}s_{1,b}s_{3,a}
\end{align*}
for the trivial reason $1\sa b=d\not=e=2\sa c$.
This gives a counterexample to the standard form of the exchange condition.


\begin{thebibliography}{Hec06b}

\bibitem[AS98]{a-AndrSchn98}
N.~Andruskiewitsch and H.-J. Schneider, \emph{Lifting of quantum linear spaces
  and pointed {H}opf algebras of order $p^3$}, J. Algebra \textbf{209} (1998),
  658--691.

\bibitem[AS05]{a-AndrSchn05p}
\bysame, \emph{On the classification of finite-dimensional pointed {H}opf
  algebras}, Preprint math.QA/0502157 (2005).

\bibitem[Bor88]{a-Borch88}
R.E. Borcherds, \emph{Generalized {K}ac-{M}oody algebras}, J. Algebra
  \textbf{115} (1988), R501--512.

\bibitem[CP61]{b-ClifPres61}
A.H. Clifford and G.B. Preston, \emph{The algebraic theory of semigroups},
  Mathematical surveys, vol.~7, Amer. Math. Soc., Providence, Rhode Island,
  1961.

\bibitem[Dri87]{inp-Drinfeld1}
V.G. Drinfel'd, \emph{Quantum groups}, Proceedings ICM 1986, Amer. Math. Soc.,
  1987, pp.~798--820.

\bibitem[Hec04]{a-Heck04e}
I.~Heckenberger, \emph{Rank 2 {N}ichols algebras with finite arithmetic root
  system}, Preprint math.QA/0412458 (2004).

\bibitem[Hec05]{a-Heck05a}
\bysame, \emph{Classification of arithmetic root systems of rank 3}, Preprint
  math.QA/0509145 (2005).

\bibitem[Hec06a]{a-Heck06b}
\bysame, \emph{Classification of arithmetic root systems}, Preprint
  math.QA/0605795 (2006).

\bibitem[Hec06b]{a-Heck04c}
\bysame, \emph{The {W}eyl groupoid of a {N}ichols algebra of diagonal type},
  Invent. Math. \textbf{164} (2006), 175--188.

\bibitem[Jim86]{a-Jimbo1}
M.~Jimbo, \emph{A $q$-analog of $\mathcal{U}(gl({N}+1))$ {H}ecke algebras, and
  the {Y}ang-{B}axter equation}, Lett. Math. Phys. \textbf{11} (1986),
  247--252.

\bibitem[Kac77]{a-Kac77}
V.G. Kac, \emph{Lie superalgebras}, Adv. Math. \textbf{26} (1977), 8--96.

\bibitem[Kac90]{b-Kac90}
\bysame, \emph{Infinite dimensional lie algebras}, Cambridge Univ. Press, 1990.

\bibitem[Kha99]{a-Khar99}
V.~Kharchenko, \emph{A quantum analog of the {P}oincar{\'e}--{B}irkhoff--{W}itt
  theorem}, Algebra and Logic \textbf{38} (1999), no.~4, 259--276.

\bibitem[LN05]{a-LoosNeh05}
O.~Loos and E.~Neher, \emph{Reflection systems and partial root systems},
  Preprint (2005).

\bibitem[Mat64]{a-Matsu64}
H.~Matsumoto, \emph{G{\'e}n{\'e}rateurs et relations des groupes de weyl
  g{\'e}n{\'e}ralis{\'e}s}, C. R. Acad. Sci. Paris \textbf{258} (1964),
  3419--3422.

\bibitem[Ser96]{a-Serg96}
V.~Serganova, \emph{On generalizations of root systems}, Commun. Algebra
  \textbf{24} (1996), 4281--4299.

\bibitem[Yam99]{a-Yam99}
H.~Yamane, \emph{On defining relations of affine {L}ie superalgebra and affine
  quantized universal enveloping superalgebras}, Publ. RIMS Kyoto Univ.
  \textbf{35} (1999), no.~3, 321--390.

\end{thebibliography}

\providecommand{\bysame}{\leavevmode\hbox to3em{\hrulefill}\thinspace}
\providecommand{\MR}{\relax\ifhmode\unskip\space\fi MR }
\providecommand{\MRhref}[2]{%
  \href{http://www.ams.org/mathscinet-getitem?mr=#1}{#2}
}
\providecommand{\href}[2]{#2}

\end{document}